\documentclass[12pt]{article}

\usepackage{amsfonts,epsfig,amsmath,amssymb}

\renewcommand{\epsilon}{\varepsilon}
\newtheorem{Theorem}{Theorem}
\newtheorem{Definition}{Definition}
\renewcommand{\phi}{\varphi}
\DeclareMathSymbol{\ophi}{\mathalpha}{letters}{"1E}

\begin{document}

\title{Stochastic Resonance in Two-State Markov
Chains}
\author{P.~Imkeller and  I.~Pavlyukevich\\ Institut
f\"ur Mathematik\\ Humboldt-Universit\"at zu
Berlin\\ Unter den Linden 6\\ 10099 Berlin\\
Germany }

\date{28 November 2000}
\maketitle

In this paper we introduce a  model
which provides a new approach to the phenomenon
of stochastic resonance. It is based on the
study of the properties of the stationary
distribution of the underlying stochastic process.
We  derive the  formula for the spectral power
amplification coefficient, study its asymptotic
properties and  dependence on parameters.


\section*{Introduction }
\label{introduction}

The notion of \emph{Stochastic Resonance} appeared
about twenty years ago in the works of Benzi et
al.~\cite{Benzi83} and  Nicolis~\cite{Nicolis82}
in the context of an attempt to explain the
phenomenon of ice ages. The modern methods of
acquiring and interpreting climate records
indicate at least seven major climate changes in
the last 700,000 years. These changes occurred with
the periodicity of about
$100,000$ years and are characterized by a
substantial variation of the average Earth's
temperature of about $10 K$.

The effect can be explained with the help of a simple energy balance model
(for an extended review on the subject
see~\cite{Imkeller00}). The Earth is considered
as a point in space, and its temporally
and spatially averaged temperature
$X(t)$ satisfies the equation
\begin{equation}
\label{en_bal_determ}
\dot X(t) = -U'(X(t)) - Q\sin{(\frac{2\,\pi
t}{T})},
\end{equation}
where $U(X)$ is a double-well potential with minima at $278.6K$ and
$288.6K$ and saddle point
at $283.3K$ and wells of equal depth.
The second term in~(\ref{en_bal_determ})
corresponds to a small variation
of the solar constant of about $0.1 \%$ with
a period of $T=100,000$ years due to the
periodic change of the eccentricity of the
Earth's orbit caused by Jupiter. The influence of
this term reflects itself in small periodic changes
of the depths of the potential wells. In this
setting, the left well is deeper during the time
intervals $(kT, (k+\frac{1}{2})T)$, whereas the
right one is deeper during the intervals
$((k+\frac{1}{2})T, (k+1)T)$, $k=0,1,2,\dots$

The trajectories of the deterministic
equation~(\ref{en_bal_determ}) have two metastable
states given by the minima of the wells. Due
to the smallness of the solar constant $Q$ no
transition between these states is possible. In
order to obtain such transitions Benzi et
al.~\cite{Benzi83} and  Nicolis~\cite{Nicolis82}
suggested to add noise to the system which results
in considering the stochastic differential equation
\begin{equation}
\label{en_bal_stoch}
\dot X^{\epsilon ,T}(t) = -U'(X^{\epsilon ,T}(t)) -
 Q\sin{(\frac{2\,\pi t}{T})}+\sqrt{\epsilon }\,\dot
W_t,
\end{equation}
$\epsilon >0$, $\dot W$ a white noise.

Now one can observe the following effect. Fix
all parameters of the system except $\epsilon $
and consider the typical behaviour of the
solutions
of~(\ref{en_bal_stoch}) for different values of
$\epsilon $. If  the noise intensity is very small,
the trajectory only occasionally can escape
from the minimum of the well in which it is
staying, and one can hardly detect any periodicity
in this motion. If the intensity is very large, the
trajectory jumps rapidly but randomly between the
two wells and therefore also lacks  periodicity
properties. An interesting effect appears when
the noise level takes a certain value
$\epsilon_0$: the trajectory always tends
to be near the minimum of the deepest well and
consequently follows the deterministic
periodic jump function which describes the
location of the deepest well's minimum. It is
very important to note that to produce this
effect one needs all three of the following
components to be present in the
system~(\ref{en_bal_stoch}): the double-well
potential for bi-stability, the noise to
pass the potential barriers, and a small
periodic perturbation to change the wells' depths.

The following are natural questions arising in the
context of these qualitative considerations: how
can one measure \emph{periodicity} of the
trajectories and, consequently, how does the
quality of \emph {tuning} of the noisy output to
the periodic input be improved by adjusting the
noise intensity
$\epsilon$?

The formulation of the latter question
suggests to consider the
system~(\ref{en_bal_stoch}) as a random amplifier.
The random system receives the
harmonic signal of small amplitude $Q$ and
usually large period
$T$ as input. The
stochastic process $X^{\epsilon ,T}(t)$ is
observed as the output. The input signal carries
power
$Q^2$ at frequency $1/T$. The random output
has continuous spectrum and thus carries
power at all frequencies.  Benzi et
al.~\cite{Benzi83} considered the power spectrum
of the output for different values of $\epsilon $
and discovered a sharp peak at the input frequency
for a certain optimal value of $\epsilon_0$. This
means that the random process
$X^{\epsilon_0,T}(t)$ has a big component
of frequency $1/T$. The effect of amplification
of the power carried by the harmonic considered as
a response of the nonlinear
system~(\ref{en_bal_stoch}) to optimally chosen
noise was called \emph{stochastic
resonance}.

In the past twenty years more than three hundred
papers on this subject were published. An
extensive  description of the phenomenon from the
physical point of view can be found in
\cite{Grammaitoni98} and
\cite{Anishchenko99}.
The notion \emph{stochastic resonance} is now used
in a much broader sense. It describes a wide class
of effects with the common underlying property: the
presence of noise induces a qualitatively new
behaviour of the system and improves some of its
characteristics.

Although stochastic resonance was observed and
studied in many physical systems, only
few  mathematically rigorous results are known. The
approach of M.~Freidlin is briefly outlined in the
next section of this paper. In sections~\ref{st_distr},
\ref{spa} and \ref{extrema} we introduce
discrete-time Markov chains with transition
probabilities chosen in such a way, that on a
large temporal scale the \emph{attractor hopping}
behaviour of the underlying diffusion process is
imitated in the limit $\epsilon\to 0$. We
investigate stochastic resonance for the Markov
chains. The last section is devoted to
generalizations and discussion.


\section{Large deviations approach}
\label{freidlin}

In this section we briefly survey rigorous
mathematical results obtained by M.~Freidlin in
\cite{Freidlin00} using the theory of large
deviations for randomly perturbed dynamical
systems, developed in Freidlin and Wentzell
(see
\cite{WentzelFreidlin84}). Though the results of
 \cite{Freidlin00} are valid in a quite general
framework, we confine our attention to a simple
example of a diffusion with weak noise.

Consider the SDE in $\mathbb R$
\begin{equation}
\label{sde}
\dot X^{\epsilon ,T}(t)= -U'(X^{\epsilon ,T}(t),
\frac{t}{T})+ \sqrt{\epsilon}\dot W(t),
\end{equation}
where $\dot W$ is a white noise and
$U'(x,t)=\frac{\partial }{\partial x}U(x,t),$
with a time dependent potential just periodically
switching between two symmetric double well
states, i.e.
\[ U(x,t)=\sum\limits_{k\geq
0}U(x){\bf 1}_{[k,k+\frac{1}{2})}(t) + U(-x){\bf
1}_{[k+\frac{1}{2},k+1)}(t),
\]
where $U(x)$ has local
minima in $x=\pm 1$  and a saddle point in $x=0$,
$\lim_{|x|\to\infty } U(x)=\infty $. We also fix
the depths of the wells by two numbers
$0<v<V,$ assuming that
$U(-1)=-V/2$,
$U(1)=-v/2$, and $U(0)=0$. Note,
that $X^{\epsilon ,T}$ is a Markov process which
is not time homogeneous. In the following Theorem
time scales are determined in which some form of
periodicity is observed.
\begin{Theorem}
Suppose $T=T(\epsilon )$ is given such that
\[
\lim\limits_{\epsilon \to 0}\epsilon \ln{T(\epsilon )}=\lambda >0.
\]

a) If $\lambda <v,$ then the Lebesgue measure of
the set
\[
\{t \in [0,1] \,:\, |X^{\epsilon
,T(\epsilon)} (T(\epsilon ) t)- \mbox{\rm
sgn}{X_0}|>\delta  \}
\]
converges to 0 in $P_{X_0}$ probability as
$\epsilon
\to 0$, for any $\delta >0$.

b) If $\lambda >v,$ then the Lebesgue measure of
the set
\[
\{t \in [0,1] \,:\, |X^{\epsilon
,T(\epsilon)} (T(\epsilon ) t)- \ophi(t)|>\delta
\}
\]
converges to 0 in $P_{X_0}$ probability as
$\epsilon
\to 0$, for any $\delta >0$, where
\[
\ophi (t)=\sum\limits_{k\geq 0}-
{\bf 1}_{[k,k+\frac{1}{2})}(t) + {\bf
1}_{[k+\frac{1}{2},k+1)}(t)
\]
and $P_{X_0}$ denotes the law of the
diffusion starting in $X_0$.
\hfill
$\square$
\end{Theorem}
It is nessesary to explain why $\lambda = v$ is
critical for the long time behaviour of the
diffusion. At least intuitively, the answer
follows from the asymptotics of the mean exit
time from a potential well for the
time-homogeneous diffusion. If the diffusion
starts in the potential well with the  depth
$v/2$, its mean time
${\bf E}(\tau(\epsilon )) $ needed to leave the
well satisfies
\[
\epsilon \ln{{\bf E}(\tau(\epsilon )) }\to v, \quad
\epsilon \to 0.
\]
according to Freidlin and Wentzell
\cite{WentzelFreidlin84}. This means, $X^{\epsilon
,T(\epsilon )}$ can leave neither the deep well
with the depth $V/2$
 nor the shallow one with the depth $v/2$ in time
$T(\epsilon )$ of order $e^{\lambda /\epsilon }$
if
$\lambda < v $. Therefore, $X^{\epsilon
,T(\epsilon )}$ stays in the $\delta
$-neighbourhood of the minimum of the  initial
well. On the other hand, if  $\lambda >v,$
$X^{\epsilon ,T(\epsilon )}$ has always enough
time to reach the deepest well. In both cases, the
Lebesgue measure of excursions leaving the  $\delta
$-tube of the deterministic periodic function
$\ophi$ is exponentially negligible on the time
scale
$T(\epsilon )$ as
$\epsilon \to 0$.

The Theorem suggests the time scale which
induces periodic and deterministic behaviour
of the system (\ref{sde}), and the  Lebesgue
measure as a measure of quality. In fact, it only
gives a lower bound for the scale. In the next
section, in the framework of discrete Markov
chains approximating the diffusion processes just
considered, we investigate different measures of
quality which provide unique optimal tuning.


\section{Markov chains with time-periodic
transition probabilities}
\label{st_distr}

For $m\in \mathbb N$, consider a Markov
chain
$X_m=(X_m(k))_{k\geq 0}$ on the state space
$\mathcal S=\{-1,1\}$.
Let $P_m(k)$ be the matrix of one-step transition
probabilities
at time $k$. If we denote $\pi^{-}_m(k)={\bf P}(X_m(k)=-1)$,
$\pi^{+}_m(k)={\bf P}(X_m(k)=1)$, and write $P^*$
for the transposed matrix, we have
\[
\left(
\begin{array}{c}
\pi^{-}_m(k+1)\\
\pi^{+}_m(k+1)
\end{array}
\right)
=P_m^*(k)
\left(
\begin{array}{c}
\pi^{-}_m(k)\\
\pi^{+}_m(k)
\end{array}
\right).
\]

In order to model the periodic switching of
the double-well potential in our Markov chains, we
define the transition matrix
$P_m$ to be periodic in time with half-period $m$.
More precisely,
\[
P_{m}(k)=
\left\{
\begin{array}{ll}
P_1, &\quad 0\leq k(\mbox{mod } 2m)\leq m-1, \\
P_2, &\quad m\leq k(\mbox{mod } 2m)\leq 2m-1,
\end{array}
\right.
\]
with
\begin{equation}
\label{P1P2}
\begin{array}{c}
P_1=
\left(
\begin{array}{cc}
1-\phi & \phi\\
\psi & 1-\psi
\end{array}
\right)
=
\left(
\begin{array}{cc}
1-px^V& px^V\\
qx^v & 1-qx^v
\end{array}
\right),
\\
P_2=
\left(
\begin{array}{cc}
1-\psi & \psi\\
\phi & 1-\phi
\end{array}
\right)
=
\left(
\begin{array}{cc}
1-qx^v& qx^v\\
px^V & 1-px^V
\end{array}
\right).
\end{array}
\end{equation}
where $\phi=\phi(\epsilon,p,V)=pe^{-V/\epsilon}$,
$\psi=\psi(\epsilon,q,v)=qe^{-v/\epsilon}$,
$x=e^{-1/\epsilon }$, $0\leq p,q\leq 1$,
$0<v<V<+\infty$, $0<\epsilon<+\infty$.  Sometimes,
it will be convenient to consider $x\in [0,1]$. In
these cases the ends of the interval will
correspond to the limits $\epsilon \to 0$ and
$\epsilon \to \infty $.

In this setting, the numbers $V/2$ and $v/2$
clearly have to be associated with the depths  of
the potential wells, $\epsilon $ with the level of
noise. According to the Freidlin-Wentzel theory,
the exponential factors in the one-step transition
probabilities just correspond to the inverses of
the expected transition times between the
respective wells for the diffusion considered in
the preceding section. This is what should be
expected for a Markov chain \emph{in equilibrium},
modulo the phenomenological \emph{pre-factors}
$p$ and
$q$. They model the pre-factors appearing in
large deviation statements, and add asymmetry to
the picture.

It is well known that for a time-homogeneous
Markov chain on $S$ with
transition matrix $P$ one can talk about 
\emph{equilibrium}, given by the stationary
distribution, to which the law of the chain
converges exponentially fast. The
stationary distribution can be found
by solving the matrix equation
 $\pi =P^*\pi $ with normalizing
condition $\pi^{-}+\pi^{+}=1$.

For non time homogeneous Markov chains
with time periodic transition matrix, the
situation is quite similar. Enlarging the
state space $S$ to $S_m = \{-1,1\}\times
\{0,1,\dots , 2m-1\},$ we recover a time
homogeneous chain by setting
\[
Y_m(k)=(X_m(k), k (\mbox{mod }2m)), \quad k\geq 0,
\]
to which the previous remarks apply. For
convenience of notation, we assume $S_m$ to be
ordered in the following way:\\ $\mathcal
S_m=[(-1,0), (1,0), (-1,1),(1,1),
\dots ,(-1,2m-1), (1, 2m-1)]$.  Writing
${\rm A}_m$ for the matrix of one-step transition
probabilities of $Y_m$, the stationary distribution
$Q=(q(i,j))^*$ is obtained as a normalized solution
of the matrix equation
$({\rm A}_m^*-E)Q=0,$ $E$ being the unit matrix.
We shall be dealing with the following variant of
stationary measure, which is not normalized in
time.
\begin{Definition}
Let $\pi _m(k)=(\pi^-_m(k),
\pi^+_m(k))^*=2m(q(-1,k),q(1,k))^*$, $0 \leq k\leq
2m-1$. We call the set $\pi_m =(\pi _m(k))_{0\leq
k\leq 2m-1}$ the stationary distribution of the
Markov chain $X_m$.
\end{Definition}
The matrix ${\rm A}_m$ of one-step transition
probabilities of $Y_m$ is
explicitly given by
\[
{\rm A}_m=\left(
\begin{array}{cccccccc}
0               &  P_1       & 0           & 0 & \cdots &0 & 0
& 0\\
0               &  0            & P_1       & 0 & \cdots &0 & 0
& 0\\
\vdots     &                 &              &    &              &   &
&\\
0               &  0            & 0           & 0 & \cdots &0 & P_2
& 0\\
0               &  0            & 0           & 0 & \cdots &0 & 0
& P_2   \\
P_2          &  0            & 0           & 0 & \cdots &0 & 0
& 0
\end{array}
\right).
\]
${\rm A}_m$ has block structure. In this notation
0 means a $2\times 2$-matrix with all entries
equal to zero, $P_1$, and $P_2$ are
the 2-dimensional matrices defined in~(\ref{P1P2}).

Applying some algebra we see that $({\rm
A}_m^*-E) Q = 0$ is equivalent to $A_m'\, Q = 0,$
where
\[
{\rm A}_m'=\left(
\begin{array}{cccccccc}
\widehat P-E  & 0       & 0      &0      & \cdots &0            & 0
& 0\\
P_1^*                             & -E      &
0      &0     & \cdots  & 0           &
0           & 0\\
\vdots                            &           &          &       &
&              &               &    \\
0                                      &  0      & 0       & 0    & \cdots
& -E         & 0            & 0\\
0                                      &  0      & 0       & 0    & \cdots &
P_2^*  & -E          & 0\\
0                                      &  0      &
0       & 0    & \cdots &0            & P_2^*   &
-E
\end{array}
\right)
\]
and $\widehat P=P_2^*P_2^*\cdots
P_1^*=(P_2^*)^m(P_1^*)^m$. But
${\rm A}_m'$ is a block-wise lower diagonal
matrix, and so $A_m' Q = 0$  can be solved in the
usual way to give
\begin{Theorem}
For every $m\geq 1$, the  stationary distribution $\pi _m$ of $X_m$ with
matrices of one-step probabilities defined in~(\ref{P1P2})  is:
\begin{equation}
\label{pi}
\begin{array}{l}
\left\{
\begin{array}{l}
\pi_m^-(l)=
\displaystyle\frac{\psi }{\phi +\psi }+
\displaystyle\frac{\phi -\psi }{\phi +\psi}\,
\frac{(1-\phi -\psi )^l}{1+(1-\phi-\psi )^m},\\
\pi_m^+(l)=
\displaystyle\frac{\phi }{\phi +\psi }-
\displaystyle\frac{\phi -\psi}{\phi +\psi}\,
\frac{(1-\phi -\psi )^l}{1+(1-\phi-\psi )^m};
\end{array}
\right.
\\
\left\{
\begin{array}{l}
\pi_m^-(l+m)=\pi _{m}^+(l),\\
\pi_{m}^+(l+m)=\pi_{m}^-(l),\qquad 0\leq l\leq m-1.
\end{array}
\right.
\end{array}
\end{equation}
\end{Theorem}
\textbf{Proof:}
$\pi _{m}(0)$ satisfies the matrix equation
$((P_2^*)^m(P_1^*)^m-E)\pi_{m}(0)=0$ with
additional condition $\pi^-_{m}(0)+\pi^+_{m}(0)=1$.
To calculate $(P_2^*)^m (P_1^*)^m$,
we use a formula for the $m$-th power of
$2\times 2$-matrices, which results in
\[
\begin{array}{l}
\left(
\begin{array}{cc}
p_{-1,-1}&p_{-1,1}\\
p_{1,-1}&p_{1,1}
\end{array}
\right)^m
{}=
\displaystyle\frac{1}{2-p_{-1,-1}-p_{1,1}}
\left(
\begin{array}{cc}
1-p_{1,1}&1-p_{-1,-1}\\
1-p_{1,1}&1-p_{-1,-1}
\end{array}
\right)\\
{}+\displaystyle
\frac{(p_{-1,-1}+p_{1,1}-1)^m}{2-p_{-1,-1}-p_{1,1}}
\left(
\begin{array}{cc}
1-p_{-1,-1}&-(1-p_{-1,-1})\\
-(1-p_{1,1})&1-p_{1,1}
\end{array}
\right)
\end{array}
\]
Using some more elementary algebra we find
\[
\begin{array}{ll}
(P_2^*)^m(P_1^*)^m&=(P_1^mP_2^m)^*=
\left(
\begin{array}{cc}
1-\psi & \psi\\
\phi & 1-\phi
\end{array}
\right)^m
\left(
\begin{array}{cc}
1-\phi & \phi\\
\psi & 1-\psi
\end{array}
\right)^m
\\
&=\displaystyle
\frac{1}{\phi+\psi}
\left(
\begin{array}{cc}
\phi & \phi\\
\psi & \psi
\end{array}
\right)
+
(1-\phi-\psi)^m\frac{\phi-\psi}{\phi+\psi}
\left(
\begin{array}{cc}
-1 & -1\\
1 & 1
\end{array}
\right)
\\
&+\displaystyle
\frac{(1-\phi-\psi)^{2m}}{\phi+\psi}
\left(
\begin{array}{cc}
\phi & -\psi\\
-\phi & \psi
\end{array}
\right),
\end{array}
\]
from which a straightforward calculation yields
\[
\left\{
\begin{array}{l}
\pi_{m}^-(0)=
\displaystyle\frac{\phi+ \psi (1-\phi -\psi )^m}{(\phi +\psi ) (1+(1-\phi
-\psi )^m)},\\
\pi_{m}^+(0)=
\displaystyle\frac{\psi +\phi (1-\phi -\psi )^m}{(\phi +\psi ) (1+(1-\phi
-\psi )^m)}.
\end{array}
\right.
\]
To compute the remaining entries, we use $\pi
_m(l)=(P_1^*)^l\pi_m(0)$ for
$0\leq l\leq m-1$,  and $\pi
_m(l)=(P_2^*)^l(P_1^*)^m\pi_m(0)$ for
$m\leq l\leq 2m-1$ to obtain~(\ref{pi}).
Note also the symmetry
$\pi_{m}^-(l+m)=\pi_{m}^+(l)$ and
$\pi_{m}^+(l+m)=\pi_{m}^-(l)$,
$0\leq l\leq m-1$.
\hfill $\square$


\section{Spectral power amplification}
\label{spa}

The chain $X_{m}$ can be interpreted as amplifier
of a signal. Our stochastic system
may be seen to receive a deterministic periodic
input signal which switches the double depths of
the potential wells in~(\ref{P1P2}), i.e.
\[
I_{m}(l)=\left\{
\begin{array}{ll}
V, &\quad 0\leq l(\mbox{mod }2m)\leq m-1,\\
v, &\quad m\leq l(\mbox{mod }2m)\leq 2m-1.
\end{array}
\right.
\]
The output is a  random process $X_{m}(k)$.

The input signal $I_{m}$ admits a spectral
representation
\[
I_{m}(k)=\frac{1}{2m}\sum_{a=0}^{2m-1}c_{m}(a)
e^{-\frac{2\pi i k}{2m}a},
\]
where
$c_{m}(a)=(1/2m)\sum_{l=0}^{2m-1}I_{2}(l)e^{\frac{2\pi
i a}{2m}l}$ is the Fourier coefficient of
frequency $a/2m$. The quantity
$|c_{2m}(a)|^2$ measures the power
carried by this Fourier component.  We are
only interested in the component of the
input frequency
$1/2m$. Its power is given by
\begin{equation}
\label{c2}
|c_{m}(1)|^2=
\frac{(V-v)^2}{4m^2}\csc^2{(\frac{\pi }{2m})}.
\end{equation}

In the stationary regime, i.e.\ if the law of $X_m$
is given by the measure $\pi_m$, the power
carried
by the output at frequency
$a/2m$ is a random variable
\[
\xi _{m}(a)=\frac{1}{2m}\sum_{l=0}^{2m-1}X_{m}(l)e^{\frac{2\pi i a}{2m}l}.
\]
We define the \emph{spectral power amplification}
as the relative expected power carried by the
component of the output with frequency
$\frac{1}{2m}$.
\begin{Definition}
The \emph{spectral power amplification
coefficient} of the Markov chain
$X_{m}$ with half period $m\geq 1$ is given by
\[
\eta_{m} =
\frac{|{\bf E}_{\pi _{m}}(\xi _{m}(1))|^2}
{|c_{m}(1)|^2}.
\]
Here
${\bf E}_{\pi _{m}}$ denotes expectation w.r.t.\ the stationary distribution
$\pi _{m}$.
\end{Definition}

The explicit description of the invariant measure
now readily yields the following formula for the
spectral power amplification.

\begin{Theorem}
Let $m\ge 1$. The spectral power amplification
coefficient of the Markov chain $X_{m}$
with  one-step  transition
probabilities~\eqref{P1P2} equals
\[
\eta_{m}=
\frac{4}{(V-v)^2}\cdot
\frac{(\phi-\psi)^2}
{(\phi+\psi)^2 +
4(1-\phi-\psi)\sin^2{(\frac{\pi}{2m})}}.
\]
\end{Theorem}
\noindent 
\textbf{Proof:}
Using \eqref{pi} one immediately gets
\begin{eqnarray*}
{\bf E}_{\pi_{m}}\xi _{m}(1) &=&
\frac{1}{2m}
\sum_{k=0}^{2m-1}{\bf
E}_{\pi_{m}}X_{m}(k)e^{\frac{2\pi i}{2m}k}=
\frac{1-e^{\pi i}}{2m}
\sum_{k=0}^{m-1}(\pi^+_{m}(k)-\pi^-_{m}(k))
e^{\frac{2\pi i }{2m}k} \\
&=& \frac{2}{m}\frac{\phi -\psi }{\phi
+
\psi }
\left(
\frac{1}{1-e^{\frac{\pi  i}{m}}}-
\frac{1}{1-(1-\phi -\psi )e^{\frac{\pi  i}{m}}}
\right).
\end{eqnarray*}
Some algebra and an appeal to \eqref{c2}
finish the proof.
\hfill $\square$

Recall now that the one-step probabilities $P_1$
and
$P_2$ depend on the parameters $0\leq p,q\leq 1$
and, what is especially important, on $0<\epsilon
<\infty $ which is interpreted as noise level. Our
next goal is to \emph{tune} the parameter $\epsilon$
to a value which maximizes the amplification
coefficient $\eta_m=\eta_m(\epsilon ) $ as a
function of $\epsilon $.


\section{Extrema and zeros of $\eta_{m}(\epsilon ) $. }
\label{extrema}

In this section we study some features of the
function $\eta_{m}(\epsilon )$ and its dependence
on $m\in \mathbb N$, $0<v<V<\infty$ and the pre-factors
$0\leq  p\leq 1$,
$0\leq q\leq 1$.

After substituting $e^{-1/\epsilon }=x$
and writing $\eta _{m}(\epsilon )=\eta _{m}(x),$
this function takes the form
\begin{eqnarray}
\label{e7}
\eta_{m} (x)=
\frac{4}{(V-v)^2}
\frac{(px^V-qx^v)^2}{(px^V+qx^v)^2 + 4(1-px^V-qx^v)\sin^2{(\frac{\pi}{2m})}}
\end{eqnarray}
In what follows, we assume $x\in [0,1]$. The
boundaries $x=0$ and $x=1$ correspond to the
limiting cases
$\epsilon = 0$ and $\epsilon = \infty $. Denote
$a_{m}=\csc{(\frac{\pi}{2m})}^2\geq 1$, $m\geq 1$.

Our main result on optimal tuning is contained in
the following theorem.
\begin{Theorem}
a) We have $\eta_{m} (x)\geq 0$, $\eta _{m}(0)=0$.

b) Let $0<\beta =\frac{v}{V}<1$ and $m\geq 1$ be fixed. There exists a continuous
function
\[
p_-(q)=p_-(q; \beta  ,m)=\displaystyle \frac{b(q; m,\beta )-
\sqrt{b(q; m, \beta )^2-4a(q; m,\beta )(2-q)q}}{2a(q;m,\beta )},
\]
where $a(q; m,\beta )= 1-a_{m}q(1-\beta )$,
$b(q;m,\beta )=2-3(1-\beta )q +a_m(1-\beta )q^2$,
with following properties:

i) $p_-(q)\geq 0$,  $q\in
[0,1]$\quad\mbox{and}\quad
$p_-(q)=0
\Leftrightarrow q=0$;

ii) $p_-(q)\leq q$, $q\in
[0,1]$\quad\mbox{and}\quad
$p_-(q)=q
\Leftrightarrow q=0$ or $q=1$, $m=1$;

iii) $\displaystyle \left. \frac{dp_-(q;m,\beta )}{dq}\right|_{q=0} =\beta $.

\noindent Moreover for $m\geq 2$

1) If $(p,q)\in U_0=\{(p,q):0< q\leq 1, 0\leq p<
p_-(q)\}$, $\eta_{m} (x)$ is strictly increasing
on $[0,1]$.

2) If $(p,q)\in U_1=\{(p,q):0< q\leq 1, p_-(q)< p\leq q\}$, $\eta_{m} (x)$
has a unique local maximum on $[0,1]$.

3) If $(p,q)\in U_2=\{(p,q):0< q\leq 1, q\leq p\leq 1\}$, $\eta_{m} (x)$
has a unique local maximum on $[0,1]$
 and a unique root on $(0,1]$. (See {\rm Fig.~})

c)For any $\delta >0$ there exists $M_0=M_0(p,q,\beta ,\delta )$ such that
for $m>M_0$
the coordinate of the local maximum $\widehat x_{m}\in[x_{m}(1-\delta ),
x_{m}]$, where
 \[
x_{m}=\left(\frac{\pi ^2}{2m^2pq}\frac{v}{V-v}\right)^{\frac{1}{V+v}}
\]
\end{Theorem}
\noindent
\textbf{Proof:}
Differentiate the explicit formula \eqref{e7}
with respect to $x$ to determine the critical
points and sets $U_0$, $ U_1$, $U_2$. The calculation
of the resonance point in $U_1$, $U_2$ requires to
find two points in some neighborhood such that
the derivative is strictly monotone on the
interval between them, and has different signs at
the extremities.
\hfill $\square$
\begin{figure}[h]
\begin{center}
\input{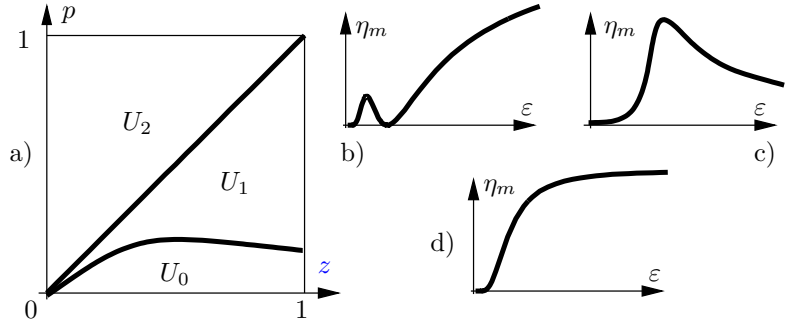}
\end{center}
\caption{a) Typical form of the domains $U_0$, $U_1$ and $U_2$.
Typical form of $\eta_m(\epsilon)$ when $(p,q)$ belongs $U_2$ (b),
$U_1$ (c) and $U_0$ (d).  }
\label{domains}
\end{figure}

\medskip
\textbf{Remarks:}\\
1. The optimal tuning rule can be rewritten in the
form
\[
m(\epsilon)\cong \frac{\pi }{\sqrt{2pq}}
\sqrt{\frac{v}{V-v}}e^{\frac{V+v}{2\epsilon }}.
\]
The maximal value of amplification is found as
\[
\lim\limits_{\epsilon \to 0}
\eta_{[m(\epsilon )]}(\epsilon ) = \frac{4}{(V-v)^2}.
\]

\noindent 2. We also see that the spectral power
amplification as a measure of quality of
stochastic resonance allows to distinguish a unique
time scale, find its exponential rate ($\lambda
=(V+v)/2$) together with the pre-exponential
factor.
\begin{figure}[h]
\begin{center}
\input{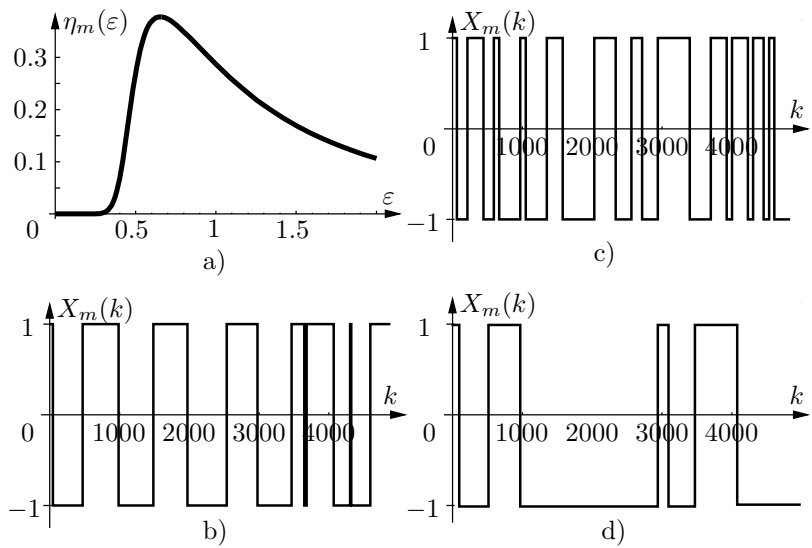}
\end{center}
\caption{a) $\eta_m(\epsilon)$ for $p=q=0.5$, $m=500$, $v=2$, $V=4$.
Numerical simulations of $X_m(k)$ for b) $\epsilon = 0.65$, 
c) $\epsilon = 0.9$ and d) $\epsilon = 0.4$. }
\label{numerics}
\end{figure}


\end{document}